\documentclass[11pt,a4paper]{article}

\usepackage{amsmath,amssymb,amscd}
\usepackage{stmaryrd}
\usepackage{graphics, psfig}
\usepackage{eufrak}

\usepackage{gastex}

\numberwithin{equation}{section}
\newcommand\refeq[1]{(\ref{#1})}

\newtheorem{deff}{D{\'e}finition}[section]

\newtheorem{prop}{Proposition}[section]
\newtheorem{property}{Property}[section]

\newenvironment{proof}{\noindent{\textsf{\underline{Proof}: }}}{$\blacksquare$}
\newtheorem{lemma}{Lemma}[subsection]
\newtheorem{theorem}{Theorem}[section]
\newtheorem{theorem*}{Theorem}

\\{}
\newenvironment{remark}{\noindent \sous{Remark}: }\\{}
\newtheorem{remark.num}{{Remark}}[section]

\newcommand\ie{i.e. }
\newcommand\eg{e.g. }

\newcommand\N{\mathbb{N}}
\newcommand\R{\mathbb{R}}
\newcommand\Z{\mathbb{Z}}

\newcommand\D{\partial}

\newcommand\dd{\text{d}}

\newcommand\lent{\llbracket }
\newcommand\rent{\rrbracket }
\newcommand\NO[1]{\ensuremath{\Arrowvert #1 \Arrowvert}}

\newcommand\vect[1]{\overrightarrow{#1}}
\newcommand\sous[1]{\underline{#1}}
\newcommand\sur[1]{\overline{#1}}
\newcommand\dsurd[2]{\frac{\partial #1}{\partial #2} }
\newcommand\chapo[1]{\widehat{#1} }
\newcommand\mc[1]{\ensuremath{\mathcal{#1}} }
\newcommand\flechegauche[1]{\overleftarrow{#1} }
\newcommand\flechedroite[1]{\overrightarrow{#1} }
\newcommand\gauchedroite[1]{\overleftrightarrow{#1} }

\title{Planar Binary Trees and Perturbative Calculus of Observables in Classical Field Theory}
\author{Dikanaina Harrivel\footnote{LAREMA, UMR 6093, Universit{\'e} d'Angers, France, \textbf{dika@tonton.univ-angers.fr}}}

\begin{document}

\maketitle

\newcommand\E{\ensuremath{E} }
\newcommand\W{\ensuremath{\mathcal{G}} }
\newcommand\lang{\left\langle }
\newcommand\rang{\right\rangle }
\newcommand\uu{\Psi }
\newcommand\RMin{\R^{n+1} }
\newcommand\DD{\sur{\Delta} }

\maketitle

\begin{abstract}
We study the Klein-Gordon equation coupled with an interaction term $(\Box+m^2)\phi+\lambda\psi^p$. 
For the linear Klein-Gordon equation, a kind of generalized Noether's theore gives us a conserved quantity. 
The purpose of this paper is to find an analogue of this conserved quantity in the interacting case. 
We see that we can do this perturbatively, and we define explicitely a conserved quantity, 
using a perturbative expansion based on Planar Binary Tree and a kind of Feynman rule. 
Only the case $p=2$ is treated but our approach can be generalized to any $\phi^p$--theory.
\end{abstract}

\section*{Introduction}
In this paper, we study the Klein--Gordon equation coupled with a second order interaction term
\begin{equation}\label{E}\tag{$E_\lambda$}
(\Box+m^2)\varphi+\lambda\varphi^2=0
\end{equation}
where $\varphi:\R^{n+1}\to\R$ is a scalar field and $\Box$ denote the operator $\frac{\D^2}{\D (x^0)^2}-\sum_{i=1}^n\frac{\D^2}{\D (x^i)^2}$. 
The constant $m$ is a positive real number which is the mass and $\lambda$ is a real parameter, the "coupling constant". Anyway our approach can be generalized 
to $\varphi^k$--theory for all $k\ge 3$. 

For any $s\in\R$ we define the hypersurface $\Sigma_s\subset\R^{n+1}$ by $\Sigma_s:=\{x=(x^0,\ldots,x^n)\in\R^{n+1}; x^0=s\}$. 
The first variable $x^0$ plays the role of time variable, and so we will denote it by $t$. Hence we interpret $\Sigma_s$ as a space--like surface 
by fixing the time to be equal to some constant $s$. 

When $\lambda$ equals zero, \refeq{E} becomes the linear Klein--Gordon equation $(\Box+m^2)\varphi=0$.
Then it is well known (see \eg \cite{ITZAK}) that for any function $\psi$ 
which satisfies $(\Box+m^2)\psi=0$, if $\varphi$ is a solution of \refeq{E} for $\lambda=0$ 
then for all $(s_1,s_2)\in\R^2$ we have
\begin{equation}\label{base}\tag{$\ast$}
\int_{\Sigma_{s_1}}\left(\dsurd{\psi}{t}\varphi-\psi\dsurd{\varphi}{t}\right)\dd\sigma
=
\int_{\Sigma_{s_2}}\left(\dsurd{\psi}{t}\varphi-\psi\dsurd{\varphi}{t}\right)\dd\sigma
\end{equation}
This last identity can be seen as expressing the coincidence on the set of solutions $\varphi$ of $(E_0)$ of two functionals $\mc{I}_{\psi,s_1}$ 
and $\mc{I}_{\psi,s_2}$
where for all function $\psi:\R^{n+1}\to\R$ and all $s\in\R$, the functional $\mc{I}_{\psi,s}$ is defined by
$$
\varphi\longmapsto\int_{\Sigma_{s}}\left(\dsurd{\psi}{t}\varphi-\psi\dsurd{\varphi}{t}\right)\dd\sigma
$$
So \refeq{base} says exactly that on the set of solutions of the linear Klein-Gordon equation, 
the functional $\mc{I}_{\psi,s}$ does not depend on the time $s$.

This could be interpreted as a consequence of a generalized version of Noether's theorem, using the fact 
that the functional 
$$
\int_K\left(\frac{1}{2}\left(\dsurd{\varphi}{t}\right)^2-\frac{\vert\nabla\varphi\vert^2}{2}-\frac{m^2}{2}\varphi^2\right)\dd x 
$$
is infinitesimally invariant under the symmetry $\varphi\to\varphi+\varepsilon\chi$, where $\chi$ is a solution of ($E_0$), up to a boundary term.

This property is no longer true when $\lambda\neq 0$ \ie when equation \refeq{E} is not linear. The purpose of this article is to obtain a
result analogous to \refeq{base} in the nonlinear (interacting) case. Another way to formulate the problem could be, if we only know the field $\varphi$ 
and its time derivative on a surface $\Sigma_{s_1}$, then how can we  evaluate $\mc{I}_{\psi,s_2}$ where $s_2\neq s_1$? 

We will see that the computation of $\mc{I}_{\psi,s_2}$ can be done perturbatively when $\lambda$ is small and $s_2$ is close to $s_1$.
This perturbative computation takes the form of a power series over Planar Binary Trees.
Note that Planar Binary Trees appear in other works on analogous Partial Differential Equations studied by perturbation (see \cite{ARBRE}, \cite{ARBRE2}, 
\cite{LEGALL}, \cite{LEGALL-DUQUESNE}, \cite{BAHNS}, \cite{Brouder.BIT}) although the point of view differs with ours.

Let us express our main result. Without loss of generality we can suppose that $s_1=0$. Then using Planar Binary Trees and starting from a function 
$\psi$ which satisfies $(\Box+m^2)\psi=0$, we explicitly construct a family of functionals $(\Psi(b))_{b\in T(2)}$ indexed by the set $T(2)$ of Planar 
Binary Trees such that the following result holds
\makeatletter
\renewcommand{\theenumi}{\roman{enumi}}
\makeatother
\begin{theorem*}\sl{
Let $q\in\N$ be such that $q>n/2$, $T>0$ be a fixed time and $\psi\in\mc{C}^2([0,T],H^{-q})$ be such that $(\Box+m^2)\psi=0$ in $H^{-q}$.
\begin{enumerate}
\item For all $\varphi$ in $\mc{C}^2([0,T],H^q)$ and $s\in[0,T]$ the power series in $\lambda$ 
\begin{equation}\label{intro}\tag{\mc{S}}
\sum_{b\in T(2)}(-\lambda)^{\vert b\vert}\lang \Psi(b)\gauchedroite{\D_s}^{\otimes \NO{b}},(\varphi,\ldots,\varphi)\rang
\end{equation}
has a non zero radius of convergence $R$. More precisely we have 
$$
R\ge 
\left(
4C_q MT
\left[
\NO{\varphi(s)}_{H^q}+\NO{\dsurd{\varphi}{t}(s)}_{H^q}
\right]
\right)^{-1}
$$
here $M$ and $C_q$ are some constants. 
\item Let $\varphi\in\mc{E}$ be such that $(\Box+m^2)\varphi+\lambda\varphi^2=0$. If the condition 
$$
8M\vert\lambda\vert C_q T\NO{\varphi}_{\mc{E}}\left(1+\vert\lambda\vert C_q T\NO{\varphi}_{\mc{E}}\right)<1
$$
is satisfied then the power series \refeq{serie} converges and we have for all $s\in[0,T]$
$$
\sum_{b\in T(2)}(-\lambda)^{\vert b\vert}\lang \Psi(b)\gauchedroite{\D_s}^{\otimes \NO{b}},(\varphi,\ldots,\varphi)\rang
=\int_{\Sigma_0}\left(\dsurd{\psi}{t}\varphi-\psi\dsurd{\varphi}{t}\right)
$$
\end{enumerate}
}
\end{theorem*}
\makeatletter
\renewcommand{\theenumi}{\arabic{enumi}}
\makeatother
The quantity $\NO{\varphi}_{\mc{E}}$ can be evaluated using the initials conditions of $\varphi$. 
More details will be available in a upcoming paper \cite{Dika2}. 
This result can be generalized for $\phi^k$--theory, $k\ge 3$, but instead of Planar Binary Trees
we have to consider Planar $k$--Trees.

Beside the fact that the functional $\sum_b \lambda^{\vert b\vert}\Psi(b)\gauchedroite{\D_s}^{\NO{b}}$ 
provides us with a kind of generalized Noether's theorem charge, it can also help us to estimate
the local values of the fields $\varphi$ and $\dsurd{\varphi}{t}$. 
We just need to 
choose the test function $\psi$ such that $\psi=0$ on the surface $\Sigma_{0}$ and $\left.\frac{\D\psi}{\D t}\right\vert_{\Sigma_{0}}$ 
is an approximation of the Dirac mass at the point $x_0\in\Sigma_{0}$. 
One gets the value 
of $\dsurd{\varphi}{t}$ at a point $x_0\in\Sigma_0$ by exchanging $\psi$ and $\dsurd{\psi}{t}$ in the previous reasoning.

Another motivation comes from the multisymplectic geometry. One of the purpose of this theory is to give a Hamiltonian formulation of the (classical) 
field theory similar to the symplectic formulation of the one dimensional Hamiltonian formalism (the Hamilton's formulation of Mechanics).
If the time variable is replaced by several space-time variables, the multisymplectic formalism is based on an analogue to the cotangent bundle, 
a manifold equipped with 
a multisymplectic form similar to the symplectic form which appears naturally in the one dimensional theory. Then starting from the Lagrangian density which 
describes the dynamics of the field, one can construct a Hamiltonian function through a Legendre transform and obtain a geometric formulation of the problem.
Note that this formalism differs from the standard Hamiltonian formulation of fields theory used by physicists (see \eg \cite{ITZAK}), in particular 
the multisymplectic approach is covariant \ie compatible with the principles of 
special and general Relativity. For an introduction to the multisymplectic geometry one can refer to \cite{sujet} and for more complete 
informations one can read the papers of F. H{\'e}lein and J. Kouneiher \cite{HK1}, \cite{HK2}.

The main motivation of the multisymplectic geometry is quantization, but it requires as preliminary to define the observable quantities, and the Poisson Bracket between 
these observables. A notion of observable have been defined by F. H{\'e}lein and J. Kouneiher. 
In the problem which interests us in this paper these observable quantities are essentially the functionals $\mc{I}_{\psi,s}$. In order to
be able to compute the Poisson bracket between two such observables $\mc{I}_{\psi_1,s_1}$ and $\mc{I}_{\psi_2,s_2}$, we must be able to 
transport $\mc{I}_{\psi_1,s_1}$ into the surface $\Sigma_{s_2}$. When $\lambda=0$, the identity \refeq{base} gives us a way to do this manipulation,
but when $\lambda\neq 0$
this is no longer the case. So F. H{\'e}lein proposed an approach based on perturbation; the reader will find more details on this 
subject in his paper \cite{sujet}.

In the first section, we begin the perturbative expansion by dealing with the linear case and the first order correction.
The second section introduces the Planar Binary Trees which
allow us to define the corrections of higher order, and the statement of the main result is given. Finally the last section contains 
the proof of the theorem.

\section{Perturbative Calculus: Beginning Expansion}
\subsection{A simple case: $\lambda=0$}\label{lambda=0}
Let us begin with the linear Klein-Gordon equation. Let $T>0$ and $s\in[0,T]$ be a fixed positive time (the negative case is similar) and $\psi:\RMin\to\R$ 
a regular function. If $\varphi$ belongs to $S_0$ \ie be a solution of the linear Klein-Gordon equation, then
\begin{equation}\label{lambda0.1}
\int_{\Sigma_{s}}\left[\dsurd{\psi}{t}\varphi-\psi\dsurd{\varphi}{t}\right]\dd\sigma-
\int_{\Sigma_{0}}\left[\dsurd{\psi}{t}\varphi-\psi\dsurd{\varphi}{t}\right]\dd\sigma
=\int_D\left[\frac{\D^2\psi}{\D t^2}\varphi-\psi\frac{\D^2\varphi}{\D t^2}\right]\dd x
\end{equation}
where $D$ denotes the set $D:=[0,s]\times\R^n\subset \R^{n+1}$. 
Since $\varphi\in S_0$ we have $\frac{\D^2\varphi}{\D t^2}=\sum_i\frac{\D^2\varphi}{\D z_i^2}-m^2\varphi$ 
hence if one replaces $\frac{\D^2\varphi}{\D t^2}$ in the right hand side of \refeq{lambda0.1} and  
perform two integrations by parts, assuming that boundary terms vanish, one obtains 
\begin{equation}\label{babar}
\int_{\Sigma_{s}}\bigl[\dsurd{\psi}{t}\varphi-\psi\dsurd{\varphi}{t}\bigr]\dd\sigma-
\int_{\Sigma_{0}}\bigl[\dsurd{\psi}{t}\varphi-\psi\dsurd{\varphi}{t}\bigr]\dd\sigma
=
\int_D\varphi(\Box +m^2)\psi
\end{equation}
Hence if we assume that $\psi$ satisfies the linear Klein-Gordon equation $(\Box +m^2)\psi=0$ 
then it follows that for all $\varphi$ in $S_0$ and for all $s\in[0,T]$  
\begin{equation}\label{I1s}
\int_{\Sigma_{s}}\left[\dsurd{\psi}{t}\varphi-\psi\dsurd{\varphi}{t}\right]\dd\sigma-
\int_{\Sigma_{0}}\left[\dsurd{\psi}{t}\varphi-\psi\dsurd{\varphi}{t}\right]\dd\sigma
=0
\end{equation}

We want to know how these computations are modified for $\varphi\in S_\lambda$ when $\lambda\neq 0$. 
For $\varphi\in S_\lambda$ we have $\frac{\D^2\varphi}{\D t^2}=\sum_i\frac{\D^2\varphi}{\D z_i^2}-m^2\varphi-\lambda\varphi^2$. 
Hence instead of \refeq{I1s} one obtains 
\begin{equation}\label{ordre1.a.virer}
\int_{\Sigma_{s}}\left[\dsurd{\psi}{t}\varphi-\psi\dsurd{\varphi}{t}\right]\dd\sigma-
\int_{\Sigma_{0}}\left[\dsurd{\psi}{t}\varphi-\psi\dsurd{\varphi}{t}\right]\dd\sigma
=\lambda\int_D\psi\varphi^2
\end{equation}
where $\psi$ is supposed to satisfy the equation $(\Box +m^2)\psi=0$. So the difference 
is no longer zero. However, one can remark that 
the difference seems\footnote{do not forget that the situation is actually more complicated because since
$\varphi$ satisfies the equation \refeq{E}, the field $\varphi$ depends on $\lambda$.} to be of order $ \lambda$. 
This is the basic observation which leads to the perturbative calculus. One can look for another functional 
which annihilates the term $\lambda\int_D\varphi^2\psi$.

\subsection{First order correction: position of the problem}
Let $s$ be a non negative integer. In the previous section, it was shown that if one choose a function $\psi$ 
such that $(\Box+m^2)\psi=0$, then equality \refeq{ordre1.a.virer} occurs for all $\varphi\in S_\lambda$. 
The purpose of this section is to search for a counter--term of order $\lambda$ which annihilates the right hand side of \refeq{ordre1.a.virer}. 

Let $\Psi^{(2)}$ be a smooth function $\Psi^{(2)}:\RMin\times\RMin\longrightarrow\R$ then for all $\varphi\in S_\lambda$ consider the quantity
\begin{equation}\label{def.I2s.prems}
\int_{\Sigma_s\times\Sigma_s}\uu^{(2)}
\left(\flechegauche{\dsurd{}{t_1}}-\flechedroite{\dsurd{}{t_1}}\right)
\left(\flechegauche{\dsurd{}{t_2}}-\flechedroite{\dsurd{}{t_2}}\right)
\varphi\otimes\varphi\ \dd \sigma_1\otimes \dd\sigma_2
\end{equation}
We need to clarify the notation $\flechegauche{A}$ and $\flechedroite{B}$ for some given operator $A$ and $B$. When the arrow 
is right to left (resp. left to right) the operator is acting on the left (resp. right).

If we assume that $\uu^{(2)}$ satisfies the boundary condition
\begin{equation}\label{condition.bords.u1}
\forall\alpha\in\{0,1\}^2\text{ , } 
\left.\frac{\D^{\vert\alpha\vert}\uu^{(2)}}{\D t^\alpha}\right\vert_{\Sigma_0\times\Sigma_s}=0
\end{equation}
then for all $\varphi$ in $S_\lambda$ we have
\begin{multline*}
\int_{\Sigma_s\times\Sigma_s}\uu^{(2)}
\left(\flechegauche{\dsurd{}{t_1}}-\flechedroite{\dsurd{}{t_1}}\right)
\left(\flechegauche{\dsurd{}{t_2}}-\flechedroite{\dsurd{}{t_2}}\right)
\varphi\otimes\varphi\\
\begin{aligned}
=&
\int_{D\times\Sigma_s}\dsurd{}{t_1}\left(
\uu^{(2)}
\left(\flechegauche{\dsurd{}{t_1}}-\flechedroite{\dsurd{}{t_1}}\right)
\left(\flechegauche{\dsurd{}{t_2}}-\flechedroite{\dsurd{}{t_2}}\right)
\varphi\otimes\varphi
\right)
\end{aligned}
\end{multline*}
here $D$ denotes the set $D:=[0,s]\times\R^n$. Assume further that we have 
\begin{equation}\label{condition.bords.u2}
\begin{array}{rcl}
\forall\alpha=(\alpha_1,\alpha_2)\in\{0,2\}\times\{0,1\}&;&\displaystyle{\left.\frac{\D^{\vert\alpha\vert}\uu^{(2)}}{\D t^\alpha}\right\vert_{D\times\Sigma_0}=0}
\end{array}
\end{equation}
then we can do the same operation for the second variable $t_2$ and finally we get
\begin{multline*}
\int_{\Sigma_s\times\Sigma_s}\uu^{(2)}
\left(\flechegauche{\dsurd{}{t_1}}-\flechedroite{\dsurd{}{t_1}}\right)
\left(\flechegauche{\dsurd{}{t_2}}-\flechedroite{\dsurd{}{t_2}}\right)
\varphi\otimes\varphi\\
=
\int_{D\times D}\uu^{(2)}\left(\flechegauche{\frac{\D^2}{\D t_1^2}}-\flechedroite{\frac{\D^2}{\D t_1^2}}\right)
\left(\flechegauche{\frac{\D^2}{\D t_2^2}}-\flechedroite{\frac{\D^2}{\D t_2^2}}\right) \varphi\otimes\varphi
\end{multline*}
Now since $\varphi$ belongs to $S_\lambda$ we have $\frac{\D^2\varphi}{\D t^2}=\sum_i\frac{\D^2\varphi}{\D z_i^2}-m^2\varphi-\lambda\varphi^2$, 
hence one can replace
the second derivatives with respect to time of $\varphi$ and then perform
integrations by parts in order to obtain : for all $\varphi$ in $S_\lambda$ the quantity \refeq{def.I2s.prems} is given by 
\begin{equation}\label{I2s.u.regulier}
\int_{D\times D}\dd x_1\dd x_2\prod_{i=1}^2\bigl(\varphi(x_i) P_i+\lambda\varphi^2(x_i)\bigr)\uu^{(2)}(x_1,x_2)
\end{equation}
where $P_i$ denotes the operator $P:=\Box+m^2$ acting on the $i$--th variable. 
Here we assume that there are no boundary terms in the integrations by parts. 

Using \refeq{I2s.u.regulier} and \refeq{ordre1.a.virer} one obtains that for all $\varphi$ in $S_\lambda$ we have 
$$
\Delta_\lambda=\lambda\biggl[\int_{D\times D}\varphi^{\otimes 2}P_1P_2\uu^{(2)}+\int_D \varphi^2\psi\biggr]+\lambda^2\cdots
$$
where $\Delta_\lambda$ denotes the quantity
\begin{multline}\label{pre.eq.u2}
\Delta_\lambda:=\int_{\Sigma_{s}}\left(\dsurd{\psi}{t}\varphi-\psi\dsurd{\varphi}{t}\right)\dd\sigma
+\lambda \int_{(\Sigma_s)^2}\uu^{(2)}
\left(\flechegauche{\dsurd{}{t_1}}-\flechedroite{\dsurd{}{t_1}}\right)
\left(\flechegauche{\dsurd{}{t_2}}-\flechedroite{\dsurd{}{t_2}}\right)
\varphi\otimes\varphi\ \dd \sigma_1\otimes \dd\sigma_2
\\-
\int_{\Sigma_{0}}\left(\dsurd{\psi}{t}\varphi-\psi\dsurd{\varphi}{t}\right)\dd\sigma
\end{multline}
Hence if we choose a function $\uu^{(2)}$ such that
\begin{equation}\label{eq.u2}
P_1P_2\uu^{(2)}(x_1,x_2)=-\delta(x_1-x_2)\psi(x_1)
\end{equation}
where $\delta$ is the Dirac operator then the first order term in the right hand side of \refeq{pre.eq.u2} vanishes. 
But because of the hyperbolicity of the operator $P$, it seems difficult to control the regularity of such a
function $\uu^{(2)}$. Hence we need to allow $\uu^{(2)}$ be in larger function space.

\subsection{Function space background}
Here we define the function spaces which allow us to express the correction terms. We fix some time $T>0$. 

Let $q\in\Z$ then we denote by $H^q(\R^n)$ (or simply $H^q$) the Sobolev space
$$
H^q(\R^n):=\left\{f\in L^2(\R^n)\ \left\vert\ (1+\vert\xi\vert^2)^{q/2}\chapo{f}(\xi)\in L^2(R^n)\right.\right\}
$$
Then it is well known (see \eg \cite{BR}, \cite{R2}, \cite{ADAMS}) that $H^q$ endowed with the norm 
$
\NO{f}_{H^q}:=\int_{\R^n} (1+\vert \xi\vert^2)^q\vert \chapo{f}\vert^2(\xi)\dd\xi
$ 
is a Banach Space. Moreover one can see in every classical text book (see \eg \cite{ADAMS}) 
the following result 
\begin{prop}\label{Hq.algebre}\sl{
If $q>n/2$ then $H^q$ is a Banach Algebra, \ie there exists some constant $C_q>0$ such that for all $(f,g)\in (H^q)^2$, $fg\in H^q$ and 
$$
\NO{fg}_{H^q}\le C_q \NO{f}_{H^q}\NO{g}_{H^q}
$$
}
\end{prop}

Until now we fix some integer $q\in\N$ such that $q>n/2$. 

\begin{deff}\sl{
Let $k$ be a positive integer, $k\in\N^*$. Then we denote by $\mc{E}^{k*}$ the space defined by 
$$
\mc{E}^{k*}:=\mc{C}^1([0,T]^k,\mc{L}_k(H^q))
$$
where $\mc{L}_k(H^q)$ denotes the space of $k$--linear continuous forms over $H^q$, we will denote $\mc{E}^{1*}$ by $\mc{E}^*$. 
}
\end{deff}
Then $\mc{E}^{k*}$ together with the norm 
$\NO{\cdot}_{k*}$ defined by
$$
\NO{U}_{k*}:=\max_{\alpha\in\{0,1\}^k}
\left(
\sup_{\substack{t\in[0,T]^k\\ (f_1,\ldots,f_k)\in (H^q)^k\\ \NO{f_j}_{H^q}\le 1}}
\left\vert
\lang 
        \frac{\D^{\vert \alpha\vert}U}{\D t^\alpha}(t),(f_1,\ldots,f_k)
\rang
\right\vert
\right)
$$ 
is a Banach Space, here $\lang\cdot,\cdot\rang$ denotes the duality brackets. 
For all $k\in\N^*$, we denote by $(\mc{E}^*)^{\otimes k}$ the space of elements $U$ of $\mc{E}^{k*}$ such that 
there exists $(U_1,\ldots,U_k)\in\mc{E}^*$ such that $U=U_1\otimes\cdots\otimes U_k$ \ie for all $(f_1,\ldots,f_k)\in (H^q)^k$ and 
for all $t=(t_1,\ldots,t_k)\in[0,T]^k$
$$
\lang U(t),(f_1,\ldots,f_k)\rang=\lang U(t_1),f_1\rang\cdots \lang U(t_k),f_k\rang
$$
Then using the fact that the space of compact supported smooth functions is dense in $H^q$, one can easily prove the following property 
\begin{property}\label{decomposable}\sl{
for all $k$ in $\N^*$,  $(\mc{E}^*)^{\otimes k}$ is a dense subspace of $\mc{E}^{k*}$.
}
\end{property}

We will denote by $\mc{E}$ the space defined by $\mc{E}:=\mc{C}^2([0,T],H^q)$. Then 
$\mc{E}$ is a Banach space and we can see naturally $\mc{E}^{*k}$ as a subspace of $\mc{L}_k(\mc{E})$ the space of $k$--linear continuous 
form over $\mc{E}$ ; $\forall U\in\mc{E}^{k*}$ and $\forall \varphi=(\varphi_1,\ldots,\varphi_k)\in\mc{E}^k$
$$
\lang U,\varphi\rang:=
\int_0^T\dd t_1\cdots\int_0^T\dd t_k \lang U(t_1,\ldots,t_k),(\varphi_1(t_1),\ldots,\varphi_k(t_k))\rang
$$
Now let us generalize the expression \refeq{def.I2s.prems} for the elements of $\mc{E}^{k*}$.
\begin{deff}\label{gauchedroite}\sl{
Let $U$ belong to $\mc{E}^*$ and $s\in[0,T]$, then we denote by $U\gauchedroite{\D_s}$ the continuous linear form over $\mc{E}$ defined by 
$\forall\varphi\in\mc{E}$
\begin{equation}\label{gd}
\lang U\gauchedroite{\D_s},\varphi\rang:=
\lang \dsurd{U}{t}(s),\varphi(s)\rang
-
\lang U(s),\dsurd{\varphi}{t}(s)
\rang
\end{equation} 
}
\end{deff}
Then using the property \ref{decomposable} one can easily prove the following property
\begin{property}\sl{
Let $k\in\N^*$ and $s\in[0,T]$ then there exists an unique operator 
$\gauchedroite{\D_s}^{\otimes k}:\mc{E}^{*k}\longrightarrow\mc{L}_k(\mc{E})$ such that for all 
$U=U_1\otimes\cdots\otimes U_k\in(\mc{E}^*)^{\otimes k}$ and for all 
$\varphi=(\varphi_1,\ldots,\varphi_k)\in\mc{E}^k$
$$
\lang \gauchedroite{\D_s}^{\otimes k}(U),\varphi\rang:=\prod_{j=1}^k\lang U_j\gauchedroite{\D_s},\varphi_j\rang
$$
For $U\in\mc{E}^{*k}$ we will denote $\gauchedroite{\D_s}^{\otimes k}(U)$ by $U\gauchedroite{\D_s}^{\otimes k}$. 
}
\end{property}

\subsection{Resolution of the first order correction}\label{section.resolution.second.ordre}
Until now we fix some function $\psi\in\mc{C}^2([0,T],H^{-q})$ such that $(\Box+m^2)\psi=0$.
In this section we will define  a functional $\Psi^{(2)}$ such that 
\begin{equation}\label{question.ordre2}
\lang\psi\gauchedroite{\D_s},\varphi\rang+\lambda\lang\Psi^{(2)}\gauchedroite{\D_s}^{\otimes 2},(\varphi,\varphi)\rang
-\lang\psi\gauchedroite{\D_0},\varphi\rang
\end{equation}
is of "order two respect with $\lambda$" for all $\varphi\in\mc{E}$ solution of \refeq{E}.

\begin{deff}\label{DD}\sl{
Let $\DD:\mc{E}^*\longrightarrow \mc{E}^{*2}$ be the operator defined by $\forall t=(t_1,t_2)\in[0,T]^2$, 
$\forall(f_1,f_2)\in (H^q)^2$
$$
\lang
\DD U(t_1,t_2),(f_1,f_2)
\rang
:=
\int_0^T\dd\tau
\lang U(\tau),\left(G(t_1)\ast f_1(\tau)\right)
\left(G(t_2)\ast f_2(\tau)\right)
\rang
$$
where for all $f\in H^q$, $t\in [0,T]$ and for all $\tau\in[0,T]$, $G(t)\ast f(\tau)$ denotes the element of $H^q$ such that 
$\forall k\in\R^n$
\begin{equation}\label{def.G}
\chapo{G(t)\ast f(\tau)}(k):=
\theta(t-\tau)\frac{\sin((t-\tau)\omega_k)}{\omega_k}\sur{\chapo{f}}(k)
\end{equation}
where $\theta$ denote the Heavyside function\footnote{$\theta(t)=0$ if $t<0$ and $1$ otherwise.} and where $\omega_k:=(m^2+\vert k\vert^2)^{1/2}$. 
}
\end{deff}
\begin{remark.num}\label{remarque.Feynman}\sl{
One can see $\DD U$ as a distribution $\DD U\in\mc{D}'((O,T)\times\R^n)$, and then one can show easily that we have the following 
formal expression for $\DD U$ 
\begin{equation}\label{DD.formel}
\DD U(x_1,x_2)=\int_{P_+}\dd y G_{ret}(x_1-y)G_{ret}(x_2-y)\psi(y)
\end{equation}
where $P_+=\{x\in\R^{n+1}\vert x^0>0\}$ and where $G_{ret}(z)$ denotes the retarded Green function of the Klein--Gordon operator
$$
G_{ret}(z):=\frac{1}{(2\pi)^n}\theta(z^0)\int_{\R^n}\dd^n k\frac{\sin(z^0\omega_{k})}{\omega_{k}}e^{ik.\sur{z}}
$$
here $\sur{z}$ denotes the spatial part of $z\in\R^{n+1}$ \ie $z=(z^0,\sur{z})$. 
}
\end{remark.num}
One can verify that $\DD$ is well defined and we have the following result 
\begin{prop}\label{prop.ordre2}\sl{
Let $\lambda$ be a real number and $s\in[0,T]$ a fixed time. 
Let $\psi\in\mc{C}^2([0,T],H^{-q})$ be such that $\frac{\D^2 \psi}{\D t^2}-\Delta\psi+m^2\psi=0$. 
If $\varphi\in\mc{E}$ is a solution of the equation \refeq{E} then the following inequality holds
\begin{multline}\label{ordre2.estimation}
\biggl\vert\lang\psi\gauchedroite{\D_s},\varphi\rang+\lambda\lang\left(\DD \psi\right)\gauchedroite{\D_s}^{\otimes 2},(\varphi,\varphi)\rang
-\lang\psi\gauchedroite{\D_0},\varphi\rang\biggr\vert\\
\le
\lambda^2\left(\frac{s^2 C_q^2}{m}\NO{\varphi}^3_{\mc{E}}+\vert\lambda\vert \frac{s^3 C_q^3}{3m^2}\NO{\varphi}^4_{\mc{E}}\right)
\NO{\psi}_{\infty,H^{-q}}
\end{multline}
}
\end{prop}
This last proposition ensures that $\Psi^{(2)}:=\DD \psi$ annihilates the term \refeq{ordre1.a.virer} of order one respect with $\lambda$. 

\begin{proof}(proposition \ref{prop.ordre2})\\
Let $\psi\in\mc{C}^2([0,T],H^{-q})$ be such that $\frac{\D^2 \psi}{\D t^2}-\Delta\psi+m^2\psi=0$ and $\varphi\in\mc{E}$ be a 
solution of equation \refeq{E}. Then since $\psi$ and $\varphi$ are $\mc{C}^2$ the function $f:t\to\lang \psi\gauchedroite{\D_t},\varphi\rang$ admits 
derivative respect with $t$ and 
$
f'(t)=\lang\frac{\D^2\psi}{\D t^2}(t),\varphi(t)\rang-\lang \psi(t),\frac{\D^2\varphi}{\D t^2}(t)\rang
$. 
But since $\psi$ and $\varphi$ satisfy $\frac{\D^2 \psi}{\D t^2}-\Delta\psi+m^2\psi=0$ and 
$\frac{\D^2\varphi}{\D t^2}-\Delta\varphi+m^2\varphi=-\lambda\varphi^2$ we have
$$
f'(t)=\lang\psi(t),(\Delta-m^2)\varphi(t)\rang-\lang \psi(t),\frac{\D^2\varphi}{\D t^2}(t)\rang
=\lambda \lang\psi(t),\varphi^2(t)\rang
$$
Hence we finally get $\forall s\in[0,T]$
\begin{equation}\label{dem.ordre2.1}
\lang\psi\gauchedroite{\D_s},\varphi\rang-\lang\psi\gauchedroite{\D_0},\varphi\rang
=f(s)-f(0)=\lambda\int_0^s \lang\psi(\tau),\varphi^2(\tau)\rang\dd\tau
\end{equation}
and we recover the identity \refeq{ordre1.a.virer} of pages \pageref{ordre1.a.virer}.

Now let us study the term of order one of the left hand side of \refeq{ordre2.estimation}. Using the definition \ref{DD} of $\DD$ 
one can show easily that it is given by the expression
\begin{equation}\label{dem.ordre2.2}
\lang\left(\DD \psi\right)\gauchedroite{\D_s}^{\otimes 2},(\varphi,\varphi)\rang
=\int_0^s\dd\tau\int_{\R^n}\dd k_1\int_{\R^n}\dd k_2
M(s,\tau,k_1)M(s,\tau,k_2)\chapo{\psi}(\tau,k_1+k_2)
\end{equation}
where $\forall (t,\tau)\in[0,T]^2$ and $\forall k\in\R^n$, the quantity $M(t,\tau,k)$ is given by
\begin{equation}\label{def.M}
M(t,\tau,k):=\cos((t-\tau)\omega_k)\sur{\chapo{\varphi}(t)}(k)-\frac{\sin((t-\tau)\omega_k)}{\omega_k}\sur{\chapo{\dsurd{\varphi(t)}{t}}}(k)
\end{equation}
The identity \refeq{dem.ordre2.2} can be seen as $\lang\left(\DD \psi\right)\gauchedroite{\D_s}^{\otimes 2},(\varphi,\varphi)\rang=u(s)$ 
where $u:[0,T]\to \R$ is the continuous function given by 
$$
u(t):=\int_0^t\dd\tau\int_{\R^n}\dd k_1\int_{\R^n}\dd k_2
M(t,\tau,k_1)M(s,\tau,k_2)\chapo{\psi}(\tau,k_1+k_2)
$$
Then in view of the definition \refeq{def.M} of $M(t,\tau,k)$ one can see that $u$ admits derivative respect with $t$ and since $u(0)=0$ we get 
$u(s)=\int_0^s u'(t)\dd t$ which leads to 
\begin{multline}\label{dem.ordre2.3}
\lang\left(\DD \psi\right)\gauchedroite{\D_s}^{\otimes 2},(\varphi,\varphi)\rang
=
\int_0^s\dd t \int_{(\R^n)^2}\dd k_1\dd k_2 \sur{\chapo{\varphi(t)}}(k_1) M(s,\tau,k_2)\chapo{\psi(t)}(k_1+k_2)\\
-
\int_0^s\dd t\int_0^t\dd\tau\int_{(\R^n)^2}\dd k_1\dd k_2  \frac{\sin((t-\tau)\omega_{k_1})}{\omega_{k_1}}\sur{\chapo{P\varphi(t)}}(k_1)
M(s,\tau,k_2)\chapo{\psi(\tau)}(k_1+k_2)
\end{multline}
where $P$ denotes the Klein--Gordon operator $P:=\Box+m^2$.

Then one can see the identity \refeq{dem.ordre2.3} as 
$\lang\left(\DD \psi\right)\gauchedroite{\D_s}^{\otimes 2},(\varphi,\varphi)\rang= v(s)+w(s)$ where 
the functions $v,w:[0,T]\to\R$ are defined by
\begin{gather*}
v(t)=\int_0^t\dd t_1 \int_{(\R^n)^2}\dd k_1\dd k_2 \sur{\chapo{\varphi(t_1)}}(k_1) M(t,\tau,k_2)\chapo{\psi(t_1)}(k_1+k_2)\\
\begin{split}
w(t)=-\int_0^t\dd t_1&\int_0^{\min(t,t_1)}\dd\tau\int_{(\R^n)^2}\dd k_1\dd k_2 \\
&\frac{\sin((t_1-\tau)\omega_{k_1})}{\omega_{k_1}}\sur{\chapo{P\varphi(t_1)}}(k_1)
M(t,\tau,k_2)\chapo{\psi(\tau)}(k_1+k_2)
\end{split}
\end{gather*}
Then one can see that $v$ and $w$ admits derivative respect with $t$ and that
\begin{multline*}
v'(t)=\int_{(\R^n)^2}\dd k_1\dd k_2
\sur{\chapo{\varphi(t)}}(k_1)\sur{\chapo{\varphi(t)}}(k_2)\chapo{\psi(t)}(k_1+k_2)\\
-\int_0^t\dd t_1\int_{(\R^n)^2}\dd k_1\dd k_2 
\frac{\sin((t-t_1)\omega_{k_2})}{\omega_{k_2}}\sur{\chapo{\varphi(t_1)}}(k_1)
\sur{\chapo{P\varphi(t)}}(k_2)\chapo{\psi(t_1)}(k_1+k_2)
\end{multline*}
and 
\begin{multline*}
w'(t)=\int_0^s\dd t_1\int_0^{\min(t_1,t)}\dd\tau \int_{(\R^n)^2}\dd k_1\dd k_2 \\
\frac{\sin((t_1-\tau)\omega_{k_1})}{\omega_{k_1}}
\frac{\sin((t-\tau)\omega_{k_2})}{\omega_{k_2}}
\sur{\chapo{P\varphi(t_1)}}(k_1)\sur{\chapo{P\varphi(t)}}(k_2)\chapo{\psi(\tau)}(k_1+k_2)\\
+\int_t^s\dd t_1\int_{(\R^n)^2}\dd k_1\dd k_2 
\frac{\sin((t_1-t)\omega_{k_1})}{\omega_{k_1}}\sur{\chapo{P\varphi(t_1)}}(k_1)
\sur{\chapo{\varphi(t)}}(k_2)
\chapo{\psi(t_1)}(k_1+k_2)
\end{multline*}
Hence since $v(0)=w(0)=0$ and using the fact that $P\varphi(t)=-\lambda\varphi^2(t)$ and in view of \refeq{def.G} 
we finally get
\begin{equation*}
\begin{split}
\biggl\langle\left(\DD \psi\right)&\gauchedroite{\D_s}^{\otimes 2},(\varphi,\varphi)\biggr\rangle=\\
&\int_0^s\dd\tau \lang \psi(\tau),\varphi^2(\tau)\rang\\
&+2\lambda\int_0^s\dd t\int_0^s\dd \tau 
\lang \psi(\tau),\left(G(t)\ast(\varphi^2(t))(\tau)\right)\varphi(\tau)\rang\\
&+\lambda^2\int_0^s\dd t_1\int_0^s\dd t_2\int_0^s\dd \tau 
\lang \psi(\tau),
\left(G(t_1)\ast(\varphi^2(t_1))(\tau)\right)
\left(G(t_2)\ast(\varphi^2(t_2))(\tau)\right)
\rang
\end{split}
\end{equation*}
Hence \refeq{dem.ordre2.1} and the last identity lead to
\begin{equation}\label{dem.ordre2.4}
\begin{split}
\lang\psi\gauchedroite{\D_s},\varphi\rang&-\lambda\lang\left(\DD \psi\right)\gauchedroite{\D_s}^{\otimes 2},(\varphi,\varphi)\rang
-\lang\psi\gauchedroite{\D_0},\varphi\rang=\\
&-2\lambda^2\int_0^s\dd t\int_0^s\dd \tau 
\lang \psi(\tau),\left(G(t)\ast(\varphi^2(t))(\tau)\right)\varphi(\tau)\rang\\
&-\lambda^3\int_0^s\dd t_1\int_0^s\dd t_2\int_0^s\dd \tau 
\lang \psi(\tau),
\left(G(t_1)\ast(\varphi^2(t_1))(\tau)\right)
\left(G(t_2)\ast(\varphi^2(t_2))(\tau)\right)
\rang
\end{split}
\end{equation}

Now to complete the proof it suffices to estimate the right hand side of \refeq{dem.ordre2.4}. Using the definition \refeq{def.G} of $G(t)\ast f$ one 
can easily prove the following lemma
\begin{lemma}\label{lemme.G}\sl{
If $f$ be in $H^q$ then for all $(t,\tau)\in[0,T]^2$ we have $(G(t)\ast f)(\tau)\in H^q$ and  
$\NO{(G(t)\ast f)(\tau)}_{H^q}\le \frac{1}{m}\theta(t-\tau)\NO{f}_{H^q}$.
}
\end{lemma}
Hence using the lemma \ref{lemme.G} and the property \ref{Hq.algebre} one get 
\begin{equation*}
\left\vert\int_0^s\dd t\int_0^s\dd \tau 
\lang \psi(\tau),\left(G(t)\ast(\varphi^2(t))(\tau)\right)\varphi(\tau)\rang\right\vert\le
\frac{s^2 C_q^2}{2m}\NO{\psi}_{\infty,H^{-q}}\NO{\varphi}_{\mc{E}}^3
\end{equation*}
and 
\begin{multline*}
\left\vert
\int_0^s\dd t_1\int_0^s\dd t_2\int_0^s\dd \tau 
\lang \psi(\tau),
\left(G(t_1)\ast(\varphi^2(t_1))(\tau)\right)
\left(G(t_2)\ast(\varphi^2(t_2))(\tau)\right)
\rang
\right\vert\le\\
\frac{s^3C_q^3}{3m^2}\NO{\psi}_{\infty,H^{-q}}\NO{\varphi}_{\mc{E}}^4
\end{multline*}
Then inserting these two inequality in \refeq{dem.ordre2.4} we finally get the result \refeq{ordre2.estimation}.
\end{proof}

Hence we found a counter--term which annihilates the term \refeq{ordre1.a.virer} of order one respect with $\lambda$. 
But some extra new terms of high order have been introduced. Thus we need to find a functional
$\lambda^2\Psi^{(3)}$ in order to delete the terms of order $\lambda^2$, and then an other functional $\lambda^3\Psi^{(4)}$ for those of order three etc. 
In order to picture all these extra terms, it will be suitable to introduce the following object:
the \emph{Planar Binary Tree}.

\section{Planar Binary Tree}
A \emph{Planar Binary Tree} (PBT) is a connected oriented tree such that each vertex has either 0 or two sons. The vertices without sons are called the \emph{leaves} and 
those with two sons are the \emph{internal vertices}. For each Planar Binary Tree, There are an unique vertex which is the son of no other vertex, this 
vertex will be called the \emph{root}. Since a Planar Binary Tree is oriented, one can define an order on the leaves. 
Until now we choose to arrange the leaves from left to right.

We will denote by $T(2)$ the set of Planar Binary Tree. 
Let denote by $\vert b\vert$ the number of internal vertices of a Planar Binary Tree $b$ and $\NO{b}$ the leave's number of $b$. Then one can easily 
show that we have $\NO{b}=\vert b\vert+1$. Let denote by $\varepsilon$ the unique Planar Binary Tree with no internal vertex. 

If $b_1$ and $b_2$ are two Planar Binary Trees, then we denote by $B_+(b_1,b_2)$ the Planar Binary Tree obtained by connecting a new root to $b_1$ on 
the left and to $b_2$ on the right. 
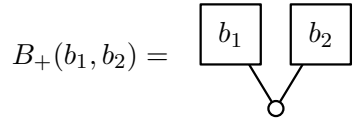
\begin{center}
\begin{picture}(0,13)(40,-20)
\gasset{linewidth=0.3,Nw=8.0,Nh=8.0,Nmr=0.0}
\node(n0)(36.11,-12.12){$b_1$}
\node(n1)(48.11,-12.12){$b_2$}
\node[linewidth=0.3,Nw=2.0,Nh=2.0,Nmr=1.0](n4)(42.11,-22){}
\gasset{AHangle=0.0}
\drawedge(n0,n4){}
\drawedge(n1,n4){}
\put(7,-16.0){$B_+(b_1,b_2)=$}
\end{picture}
\end{center}
Then one can easily show that $\vert B_+(b_1,b_2)\vert=\vert b_1\vert+\vert b_2\vert +1$ and $\NO{B_+(b_1,b_2)}=\NO{b_1}+\NO{b_2}$, and 
for all $b\in T(2)$, $b\neq\varepsilon$, there is an unique couple $(b_1,b_2)\in T(2)^2$ such that 
$b=B_+(b_1,b_2)$. For further details on the Planar Binary Trees, one can consult \cite{GKP}, \cite{SEDGEWICK},
\cite{FRA-simplicial} or \cite{LAAN}. 

\begin{deff}\label{def.Delta}\sl{
We define inductively the family $(\DD(b))_{b\in T(2)}$ of functionals $\DD(b):\mc{E}^*\longrightarrow\mc{E}^{*\NO{b}}$ by
\begin{equation}\label{deff.DDb}
\left\{
\begin{array}{l}
\displaystyle{\DD(\varepsilon):=id}\\
\displaystyle{\forall (b_1,b_2)\in T(2)^2\text{ ; }\DD(B_+(b_1,b_2)):=(\DD(b_1)\otimes\DD(b_2))\circ\DD}
\end{array}
\right.
\end{equation}
where for $\mc{U}:\mc{E}^*\to\mc{E}^{*k}$ and $\mc{V}:\mc{E}^*\to\mc{E}^{*l}$, $\mc{U}\otimes \mc{V}$ denotes the unique functional 
from $\mc{E}^{*2}$ to $\mc{E}^{*(k+l)}$ such that for all $U=U_1\otimes U_2\in (\mc{E}^*)^{\otimes 2}$, 
$\mc{U}\otimes \mc{V}(U)=\mc{U}(U_1)\otimes \mc{V}(U_2)$. 
}
\end{deff}
Let $\psi$ belong to $\mc{E}^*$, then we consider the family $(\Psi(b))_{b\in T(2)}$ defined by 
$$
\Psi(b):=\DD(b)(\psi)\in \mc{E}^{*\NO{b}}
$$

Then using remark \ref{remarque.Feynman}  we can see that formally the functionals $\Psi(b)$, $b\in T(2)$, can be constructed using the 
following rules :
\begin{enumerate}
\item attach to each leaf of $b$ the space--time variable $x_1$, $x_2$, $\ldots,x_{\NO{b}}$  
with respect to the order of the leaves.
\item for each internal vertex attach a space--time integration variable $y_i\in\R^{n+1}$ and 
integrate this variable over $P_+$.
\item for each line between the vertices $v$ and $w$ where the depth of $v$ is lower than the $w$'s, put a factor $G_{ret}(a_v-a_w)$ 
where $a_v$ (resp. $a_w$) is the space--time variable associated with $v$ (resp. $w$).
\item finally multiply by $\psi(a_r)$ where $a_r$ is the space--time variable attached to the root of the Planar Binary Tree $b$.
\end{enumerate}

To fix the ideas, let us treat an example. Let
$b\in T_3$ be the Planar Binary Tree described by the following graph
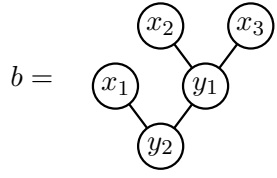
\begin{center}
\begin{picture}(22,18)(0,-28)
\put(-10,-20.0){$b =$}
\gasset{linewidth=0.3,Nw=6.0,Nh=6.0,Nmr=3.0,AHangle=0.0,AHLength=1.55}
\node(n3)(4,-20){$x_1$}
\node(n4)(10.0,-28.0){$y_2$}
\node(n5)(16.0,-20.0){$y_1$}
\node(n6)(10.0,-12.0){$x_2$}
\node(n7)(22.0,-12.0){$x_3$}
\drawedge(n4,n3){}
\drawedge(n4,n5){}
\drawedge(n5,n6){}
\drawedge(n5,n7){}
\end{picture}
\end{center}
Then using definition \ref{def.Delta} we have $\Psi(b)=(id\otimes\DD)\circ\DD\psi$ and formally for 
$x=(x_1,x_2,x_3)\in ([0,T]\times\R^n)^3$, $\Psi(b)(x)$ is given by the following 
\begin{equation*}
\Psi(b)(x)=\iint_{P_+}\dd y_1\dd y_2
G_{ret}(x_1-y_2)G_{ret}(y_1-y_2)G_{ret}(x_2-y_1)G_{ret}(x_3-y_1)\psi(y_2)
\end{equation*}

\makeatletter
\renewcommand{\theenumi}{\roman{enumi}}
\makeatother
\begin{theorem}\label{prop.fin}\sl{
\begin{enumerate}
\item Let $\psi\in\mc{C}^2([0,T],H^{-q})$ be such that $(\Box+m^2)\psi=0$ in $H^{-q}$. Let $\varphi$ be in $\mc{E}$ and $s\in[0,T]$, 
then the power series in $\lambda$ 
\begin{equation}\label{serie}\tag{$\ast$}
\sum_{b\in T(2)}(-\lambda)^{\vert b\vert}\lang \Psi(b)\gauchedroite{\D_s}^{\otimes \NO{b}},(\varphi,\ldots,\varphi)\rang
\end{equation}
has a non zero radius of convergence $R$. More precisely we have 
$$
R\ge 
\left(
4C_q MT
\left[
\NO{\varphi(s)}_{H^q}+\NO{\dsurd{\varphi}{t}(s)}_{H^q}
\right]
\right)^{-1}
$$
here $M$ is defined by $M:=\max(\frac{1}{m},1)$ and $C_q$ is the constant of the property \ref{Hq.algebre}.
\label{prop.a}
\item Let $\varphi\in\mc{E}$ be such that $(\Box+m^2)\varphi+\lambda\varphi^2=0$. If the condition 
\begin{equation}\label{condition.fin}
8M\vert\lambda\vert C_q T\NO{\varphi}_{\mc{E}}\left(1+\vert\lambda\vert C_q T\NO{\varphi}_{\mc{E}}\right)<1
\end{equation}
is satisfied then the power series \refeq{serie} converges and we have for all $s\in[0,T]$
$$
\sum_{b\in T(2)}(-\lambda)^{\vert b\vert}\lang \Psi(b)\gauchedroite{\D_s}^{\otimes \NO{b}},(\varphi,\ldots,\varphi)\rang
=\lang \psi\gauchedroite{\D_0},\varphi\rang
$$
\label{prop.b}
\end{enumerate}
}
\end{theorem}
\makeatletter
\renewcommand{\theenumi}{\arabic{enumi}}
\makeatother

\begin{remark}
Note that it is possible to control the norm $\NO{\varphi}_{\mc{E}}$ with the norm of initials data. More preciselly for all 
any $(\varphi^0,\varphi^1)\in (H^q)^2$, $\lambda\in\R$ and $T\in\R$ such that $T\vert \lambda\vert \NO{(\varphi^0,\varphi^1)}$ is 
small enough, it is possible to construct a solution $\varphi\in\mc{C}^2([0,T],H^q)$ of \refeq{E} such that 
 $\varphi(0,\cdot)=\varphi^0$ and $\dsurd{\varphi}{t}(0,\cdot)=\varphi^1$. Then one can control $\NO{\varphi}_{\mc{E}}$ using $\NO{(\varphi^0,\varphi^1)}$. 
A proof of this result, based on a remark of Christian Brouder (\cite{Brouder.BIT}) will be expand in \cite{Dika2}.
\end{remark}

Let us comment this last proposition. First of all, using definition \refeq{gauchedroite} of $\gauchedroite{\D_s}$, 
one can remark that the power series \refeq{serie} depends only on 
$\varphi(s,\cdot)$ and $\dsurd{\varphi}{t}(s,\cdot)$. 
Hence the theorem answers the original question. 

We have written the solution for $s$ non negative, but the study can be done in the same way for negative $s$. 
Finally the result exposed in the proposition \ref{prop.fin} can be generalized to $\phi^p$--theory \ie for 
the equation
$(\Box+m^2)\varphi+\lambda\varphi^p=0$, $p\ge 2$. But the set of Planar Binary Trees must be replaced by T(p), the set of 
\emph{Planar $p$-Trees} \ie oriented rooted trees which vertices have 
$0$ or $p$ sons, then the definition of $(\Psi(b))_{b\in T(p)}$ remains the same \ie $\Psi(b):=\DD^{(p)}(b)\psi$ where $\DD^{(p)}(b)$ is an adaptation 
of definition \ref{def.Delta} for $p$--trees. Then an analogue of theorem \ref{prop.fin} holds but the condition \refeq{condition.fin} must be adapted. 

\section{Proof of the main proposition}
\subsection{Radius of convergence}\label{premiere.partie}
Let us deal with the first part of theorem \ref{prop.fin}. 
First we will prove the following lemma
\begin{lemma}\label{lemme.RdC}\sl{
Let $\psi$ belong to $\mc{E}^*$ and $b\in T(2)$, then we have  
\begin{equation}\label{base.dem.RdC}
\NO{\DD(b)\psi}_{*\NO{b}}\le \left( C_q M T\right)^{\vert b\vert}\NO{\psi}_{*1}
\end{equation}
where $M$ and $C_q$ are the constants which appear in theorem \refeq{prop.fin}.
}
\end{lemma}
\begin{proof}(lemma \ref{lemme.RdC})\\
We will show \refeq{base.dem.RdC} inductively with respect to $\vert b\vert$ the number of internal vertices of $b$.

If $b=\varepsilon$ then inequality \refeq{base.dem.RdC} is satisfied. Let $N\in\N$, suppose that \refeq{base.dem.RdC} is true 
for all $b\in T(2)$ such that $\vert b\vert\le N$ and let $b\in T(2)$ be such that $\vert b\vert=N+1\ge 1$. Then $b$ writes 
$b=B_+(b_1,b_2)$ and by definition we have $\Psi(b)=\DD(b)\psi=(\DD(b_1)\otimes\DD(b_2))\circ\DD\psi$. 
Since $(\mc{E}^*)^{\otimes 2}$ is dense in $\mc{E}^{*2}$ there is a sequence 
$U_n=U^{(1)}_n\otimes U^{(2)}_n\in (\mc{E}^*)^{\otimes 2}$, $n\in\N$, such that $U_n\to \DD\psi$ in $\mc{E}^{*2}$. Then 
one can show easily that 
$$
\NO{\DD(b_1)\otimes\DD(b_2)) U_n}_{*(\NO{B_+(b_1,b_2)})}=\NO{\DD(b_1)U_n^{(1)}}_{\NO{b_1}}\NO{\DD(b_2)U_n^{(2)}}
$$
but since $\vert B_+(b_1,b_2)\vert=\vert b_1\vert+\vert b_2\vert+1$ we have $\vert b_1\vert\le N$ and 
$\vert b_2\vert\le N$, hence \refeq{base.dem.RdC} is valid for $b_1$ and $b_2$ we finally get 
\begin{equation*}
\begin{split}
\NO{\DD(b_1)\otimes\DD(b_2)) U_n}_{*(\NO{B_+(b_1,b_2)})}
\le& 
\left( C_q M T\right)^{\vert b_1\vert+\vert b_2\vert}\NO{U^{(1)}_n}_{*1}\NO{U^{(2)}_n}_{*1}\\
=&
\left( C_q M T\right)^{\vert b_1\vert+\vert b_2\vert}\NO{U^{(1)}_n\otimes U^{(2)}_n}_{*2}
\end{split}
\end{equation*}
Then taking the limit $n\to\infty$ in the previous inequality leads to 
\begin{equation}\label{dem.rec.RdC}
\NO{\DD(B_+(b_1,b_2))\psi}_{*(\NO{B_+(b_1,b_2)})}
\le 
\left( C_q M T\right)^{\vert b_1\vert+\vert b_2\vert}
\NO{\DD\psi}_{2*}
\end{equation}

Let $f_1$ and $f_2$ belong to $H^q$ then by definition we have for all $(t_1,t_2)\in[0,T]^2$ and $\alpha=(\alpha_1,\alpha_2)\in\{0,1\}^2$
\begin{equation}\label{froufrou}
\lang \frac{\D^{\vert\alpha\vert}\DD\psi}{\DD t^\alpha}(t_1,t_2),(f_1,f_2)\rang
=
\int_0^T\dd \tau 
\lang \psi(\tau),
\left(G^{\alpha_1}(t_1)f_1(\tau)\right)\left(G^{\alpha_2}(t_2)f_2(\tau)\right)\rang
\end{equation}
where $G^0(t)f(\tau):=G(t)f(\tau)$ has been defined in the section \ref{section.resolution.second.ordre} page \pageref{def.G} and where 
$G^1(t)f(\tau)\in H^q$ is the function such that $\forall \vect{k}\in \R^n$
$$
\chapo{G^1(t_1)f(\tau)}(\vect{k}):=\theta(t_1-\tau)\cos((t_1-\tau)\omega_{\vect{k}})\chapo{f}(k)
$$
Then one can easily show that $\NO{G^1(t)f(\tau)}_{H^q}\le \frac{1}{m}\theta(t-\tau)\NO{f}_{H^q}$ and 
$\NO{G^1(t)f(\tau)}_{H^q}\le \theta(t-\tau)\NO{f}_{H^q}$. Hence inserting these results in \refeq{froufrou} and using 
property \ref{Hq.algebre} we get 
$$
\NO{\DD\psi}_{2*}\le MC_q T \NO{\psi}_{*1}
$$
So in view of \refeq{dem.rec.RdC} we see that the estimation \refeq{base.dem.RdC} is valid for $b=B_+(b_1,b_2)$. 
\end{proof}

Let $\varphi$ belong to $\mc{E}$ then lemma \ref{lemme.RdC} shows that for all $s\in[0,T]$ and for all $b\in T(2)$ we have 
$$
\left\vert
\lang\Psi(b)\gauchedroite{\D_s},(\varphi,\ldots,\varphi)\rang
\right\vert
\le
\left(C_q MT\right)^{\vert b\vert}\NO{\psi}_{*1}
\left[
\NO{\varphi(s)}_{H^q}+\NO{\dsurd{\varphi}{t}(s)}_{H^q}
\right]^{\NO{b}}
$$
then using the fact (see \eg \cite{SEDGEWICK}) that the number $p_N$ of Planar Binary Tree $b$ such that $\vert b\vert=N$ satisfies 
$p_N\le 4^N$ we finally get the first part of theorem \ref{prop.fin}, \ie the power series in $\lambda$ defined by 
$\sum_{b\in T(2)}(-\lambda)^{\vert b\vert }\lang\Psi(b)\gauchedroite{\D_s},(\varphi,\ldots,\varphi)\rang$ has a non--zero 
radius of convergence $R$ and 
$$
R\ge 
\left(
4C_q MT
\left[
\NO{\varphi(s)}_{H^q}+\NO{\dsurd{\varphi}{t}(s)}_{H^q}
\right]
\right)^{-1}>0
$$

\subsection{Algebraic calculations}
Let us fix some time $s$ in $[0,T]$, then we define the operator $P:\mc{E}^*\longrightarrow\mc{F}'$ where $\mc{F}\subset\mc{E}$ denotes the space 
$
\mc{F}:=\mc{C}^2([0,T],H^{q})\cap\mc{C}^0([0,T],H^{q+2})
$ by 
for all $U\in\mc{E}^*$ and for all $\varphi\in\mc{F}$
\begin{equation}\label{def.P}
\lang PU,\varphi\rang:=\lang U\gauchedroite{\D_s},\varphi\rang-\lang U\gauchedroite{\D_0},\varphi\rang
+\int\dd\tau \lang U(\tau),(\Box+m^2)\varphi(\tau)\rang
\end{equation}
here $\Box+m^2$ denotes the operator $\mc{F}\to\mc{E}$ defined by $\Box=\frac{\D^2}{\D t^2}-\Delta$. Let $k$ be an integer 
$k\in\N^2$ then for all $I\subset\lent 1,k\rent$ we denote by $P^k_I$ the unique continuous operator 
$P^k_I:\mc{E}^{k*}\longrightarrow\mc{L}_k(\mc{F})$ such that for all $U=U_1\otimes\cdots U_k\in(\mc{E}^*)^{\otimes k}$ and 
for all $\varphi=(\varphi_1,\ldots,\varphi_k)\in\mc{F}^k$
$$
\lang P^k_I U,\varphi\rang=
\prod_{i\in I}\lang P U_i,\varphi_i\rang \prod_{j\not\in I}\int_0^T \lang U_j(\tau_j),\varphi_j(\tau_j)\rang\dd\tau_j
$$

Let $\varphi\in\mc{E}$ be a solution of \refeq{E} \ie such that $(\Box+m^2)\varphi=-\lambda\varphi^2$. Then in view property \ref{Hq.algebre} 
$\varphi$ belongs to $\mc{F}$. Let $b$ be a Planar Binary Tree such that $b\neq\varepsilon$ and denote by $k$ the leave's number of $b$, 
$k:=\NO{b}$. Then in view of definition \ref{def.Delta} one can easily 
see that for all $J\subset\lent 1,k\rent$ we have $\Psi(b)\gauchedroite{\D_0^I}=0$ hence the definition \refeq{def.P} of $P$ leads to 
\begin{equation}\label{base.dem0}
\lang \Psi(b)\gauchedroite{\D_s}^{\otimes k},(\varphi,\ldots,\varphi)\rang
=\sum_{I\subset\lent 1,k\rent}
\lambda^{k-\vert I\vert}\lang P^{k}_I \Psi(b),(\varphi^{\alpha^I_1},\ldots,\varphi^{\alpha^I_k})\rang
\end{equation}
where $\alpha^I_j=2$ if $j\not \in I$ and $\alpha^I_j=1$ otherwise;  
here we use the fact that $\varphi$ satisfies $-(\Box+m^2)\varphi=\lambda\varphi^2$. Moreover the proof of proposition \ref{prop.ordre2} shows 
that if one choose $\psi\in\mc{C}^2([0,T],H^{-q})$ such that $(\Box+m^2)\psi=0$ then 
we have 
\begin{equation}\label{base.dem0.ordre0}
\lang \psi\gauchedroite{\D_s},\varphi\rang-
\lang \psi\gauchedroite{\D_0},\varphi\rang
=
-\lambda \lang \psi,\varphi^2\rang=-\lambda \int_0^s\lang \psi(\tau),\varphi^2(\tau)\rang\dd\tau
\end{equation}
\ie $P\psi=0$. For $N\in\N^*$ let denote by $\Delta_N$ the finite sum 
$$
\Delta_N:=\sum_{\substack{b\in T(2)\\ \vert b\vert\le N}}(-\lambda)^{\vert b\vert}\lang \Psi(b)\gauchedroite{\D_s}^{\otimes \NO{b}},(\varphi,\ldots,\varphi)\rang
-\lang \psi\gauchedroite{\D_0},\varphi\rang
$$
Then \refeq{base.dem0.ordre0} and \refeq{base.dem0} lead to 
\begin{equation}\label{base.dem}
\Delta_N=
\sum_{\beta=1}^{2N-1}\lambda^{\beta-1}
\sum_{\substack{1\le k\le N\\ 0\le l\le k\\ k+l=\beta}}
\sum_{\substack{b\in T(2)\\ \NO{b}=k}}
\sum_{\substack{I\subset \lent 1,k\rent\\ \vert I\vert =k-l}}
(-1)^{\vert b\vert}
\lang P^k_I \Psi(b),(\varphi^{\alpha^I_1},\ldots,\varphi^{\alpha^I_k})\rang
\end{equation}

Let $\beta\in\N^*$ be such that $\beta\le N$ then $\Delta_N^\beta$ the term of order $\beta$ with respect to $\lambda$ in \refeq{base.dem} writes
\begin{multline}\label{base.dem.beta<N}
\Delta_N^\beta=
\sum_{\substack{b\in T(2)\\ \NO{b}=\beta}}
(-1)^{\vert b\vert}\lang P^\beta_{\lent 1,\beta\rent}\Psi(b),(\varphi,\ldots,\varphi)\rang
\\
+
\sum_{\substack{1\le l\le k\le \beta\\ k+l=\beta}}
\sum_{\substack{a\in T(2)\\ \NO{a}=k}}
\sum_{\substack{I\subset \lent 1,k\rent\\ \vert I\vert =k-l}}
(-1)^{\vert a\vert}
\lang P^k_I \Psi(a),(\varphi^{\alpha^I_1},\ldots,\varphi^{\alpha^I_k})\rang
\end{multline}
Let us focus on the first sum of this last identity. 
We need some extra structure on the set of Planar Binary Tree. 
Two special Planar Binary Trees play an important role : the Planar Binary Tree $\varepsilon$ with one leaf, and $Y$ the one with two leaves
\begin{center}
\begin{picture}(67,4)(0,-10)
\gasset{Nw=2.0,Nh=2.0,Nmr=1.0}
\put(1.0,-8.5){$\varepsilon=$}
\node(n1)(12.0,-8.0){}
\put(32.0,-8.0){and}
\put(52.0,-8.0){$Y=$}
\node(n4)(61.0,-5.0){}
\node(n5)(64.0,-11.0){}
\node(n6)(67.0,-5.0){}
\gasset{AHangle=0.0,AHLength=1.65}
\drawedge(n5,n4){}
\drawedge(n6,n5){}
\end{picture}
\end{center}
We also define the \emph{growing} operation. Let $b$ be a Planar Binary Tree with $k$ leaves and $E=(E_1,\ldots,E_k)$ be a 
$k$-uplet in $\{\varepsilon,Y\}^k$. 
We call the \emph{growing of $E$ on $b$} and denote by $E\propto b$ the Planar Binary Tree obtained by replacing the $i$-th leaf of $b$ by
$E_i$. As an example
\begin{center}
\begin{picture}(137,12)(0,-14)
\gasset{Nw=2.0,Nh=2.0,Nmr=1.0,AHangle=0.0,AHLength=1.6}
\put(5.0,-8.0){$\bigl($}
\node(n1)(9.0,-5.0){}
\node(n2)(12.0,-11.0){}
\node(n3)(15.0,-5.0){}
\drawedge(n1,n2){}
\drawedge(n2,n3){}
\put(20.0,-8.0){,}
\node(n4)(23,-5.0){}
\node(n5)(26.0,-11.0){}
\node(n6)(29.0,-5.0){}
\drawedge(n4,n5){}
\drawedge(n6,n5){}
\put(34.0,-8.0){$\bigr)$}

\put(39.0,-8.0){$\propto$}

\node(n7)(44,-5.0){}
\node(n8)(47.0,-11.0){}
\node(n9)(50.0,-5.0){}
\drawedge(n7,n8){}
\drawedge(n8,n9){}

\put(57.0,-8.0){$=$}

\put(62.0,-8.0){$\bigl($}
\node(n9)(66.0,-8.0){}
\put(68.0,-8.0){$,$}
\node(n10)(72.0,-8.0){}
\put(74.0,-8.0){,}
\node(n11)(77.0,-5.0){}
\node(n12)(80.0,-11.0){}
\node(n13)(83.0,-5.0){}
\drawedge(n11,n12){}
\drawedge(n12,n13){}
\put(88.0,-8.0){$\bigr)$}

\put(93.0,-8.0){$\propto$}

\node(n25)(100,-8.0){}
\node(n26)(103,-14){}
\node(n27)(106,-8.0){}
\node(n30)(97,-2){}
\node(n31)(103,-2){}
\drawedge(n30,n25){}
\drawedge(n31,n25){}
\drawedge(n26,n25){}
\drawedge(n26,n27){}

\put(114.0,-8.0){$=$}

\node(n13)(121,-8.0){}
\node(n14)(126,-14){}
\node(n16)(118,-2){}
\node(n17)(124,-2){}
\drawedge(n16,n13){}
\drawedge(n17,n13){}
\drawedge(n14,n13){}
\node(n18)(128.0,-2){}
\node(n19)(131.0,-8.0){}
\node(n20)(134.0,-2){}
\drawedge(n18,n19){}
\drawedge(n19,n20){}
\drawedge(n14,n19){}

\end{picture}
\end{center}
For $E\in\{\varepsilon,Y\}^k$ we denote by $n_Y(E)$ the occurrence number of $Y$ in $E$  \ie $n_Y(E):=\text{Card}\{i\vert E_i=Y\}$. Then we have the
 combinatorial lemma 
\begin{lemma}\label{noeud.dem}\sl{
\begin{enumerate}
\item
Let $b$ be a Planar Binary Tree with $\beta$ leaves, $\beta\ge 2$. Then we have
$$
-\sum_{\substack{1\le l\le k\le \beta\\k+l=\beta}}\sum_{\substack{a\in T_k\\E\in\{\varepsilon,Y\}^k\vert  n_Y(E)=l\\\text{such that }E\propto a=b}}
(-1)^{\vert a\vert}=(-1)^{\vert b\vert}
$$
\label{point1.lemme}
\item Let $p\in\N^*$, $a\in T(2)$ be such that $\NO{a}=p$ and $E\in\{\varepsilon, Y\}^p$ then we have $\NO{E\propto a}=p+n_Y(E)$ and
$$
\lang P^{p+n_Y(E)}_{\lent 1,p+n_Y(E)\rent}\Psi(E\propto a), (\varphi,\ldots,\varphi)\rang=
\lang P^p_{I_E}\Psi(a),(\varphi^{\alpha^{I_E}_1},\ldots,\varphi^{\alpha^{I_E}_p}
\rang
$$
where $I_E:=\{j\in\lent 1,p\rent\text{ such that }E_j=\varepsilon\}$ and $\alpha^{ I_E}_j=2$ if $j\not\in I_E$ and $1$ otherwise. 
\label{point2.lemme}
\end{enumerate}
}
\end{lemma}
We postpone the proof until appendix. Then the point \ref{point1.lemme} lemma \ref{noeud.dem} leads to
\begin{multline}\label{root.dem}
\sum_{\substack{b\in T(2)\\ \NO{b}=\beta}}(-1)^{\vert b\vert}\lang P^\beta_{\lent 1,\beta\rent}\Psi(b),(\varphi,\ldots,\varphi)\rang\\
=-\sum_{\substack{1\le l\le k\le \beta\\k+l=\beta}}\sum_{\substack{a\in T_k\\E\in\{\varepsilon,Y\}^k\vert  n_Y(E)=l}}
(-1)^{\vert a\vert}
\lang P^\beta_{\lent 1,\beta\rent}\Psi(E\propto a),(\varphi,\ldots,\varphi)\rang
\end{multline}
But since $E\in\{\varepsilon, Y\}^p$ is entirely determined by $p$ and $I_E$, the point \ref{point2.lemme} of 
lemma \ref{noeud.dem} and identity \refeq{root.dem} lead to 
\begin{multline*}
\sum_{\substack{b\in T(2)\\ \NO{b}=\beta}}(-1)^{\vert b\vert}\lang P^\beta_{\lent 1,\beta\rent}\Psi(b),(\varphi,\ldots,\varphi)\rang\\
=-\sum_{\substack{1\le l\le k\le \beta\\k+l=\beta}}
\sum_{\substack{a\in T(2)\\ \NO{a}=k}}
\sum_{\substack{I\subset \lent 1,k\rent\\ \vert I\vert=k-l}}
(-1)^{\vert a\vert}
\lang P^k_I\Psi(a),(\varphi^{\alpha^I_1},\ldots,\varphi^{\alpha^I_k})\rang
\end{multline*}
then inserting this last identity in \refeq{base.dem.beta<N} we finally get that for all $\beta\le N$, $\Delta_\beta^N=0$. 

Hence we have shown that $\Delta_N=\sum_{\beta=1}^{2N-1}\lambda^{\beta-1}\Delta_\beta^N$ is of order $N$ with respect to $\lambda$. To complete 
the proof it suffices to show that $\Delta_N$ converges to $0$ when $N$ tends to infinity.

\subsection{analytic study}
We have shown that all the terms of order $\beta$ with $\beta\le N$ in identity \refeq{base.dem} vanish, hence we have
\begin{equation}\label{base.dem.2}
\Delta_N=
\sum_{\beta=N+1}^{2N-1}\lambda^{\beta-1}
\sum_{\substack{1\le l\le k\le N\\ k+l=\beta}}
\sum_{\substack{b\in T(2)\\ \NO{b}=k}}
\sum_{\substack{I\subset \lent 1,k\rent\\ \vert I\vert =k-l}}
(-1)^{\vert b\vert}
\lang P^k_I \Psi(b),(\varphi^{\alpha^I_1},\ldots,\varphi^{\alpha^I_k})\rang
\end{equation}
We have to estimate the right hand side of this last identity. Let us prove the following lemma
\begin{lemma}\label{lemme.maj.PPsi}\sl{
Let $k\in\N^*$, $k\ge 2$ and $b\in T(2)$ be such that $\NO{b}=k$, then for all $I\subset\lent 1,k\rent$, $\varphi\in\mc{E}$ solution of \refeq{E} and 
$\psi\in\mc{E}^*$ we have
\begin{multline}\label{maj.fin}
\left\vert
\lang P^k_I\Psi(b),(\varphi^{\alpha^I_1},\ldots,\varphi^{\alpha^I_k})\rang
\right\vert\le\\
(C_qT)^{k-\vert I\vert}
\NO{\varphi}^{2k-\vert I\vert}_{\mc{E}}(C_qMT)^{k-1}
\left(2+\vert \lambda\vert C_qT\NO{\varphi}_{\mc{E}}\right)^{\vert I\vert}\NO{\psi}_{\mc{E}^*}
\end{multline}
}
\end{lemma}
\begin{proof}(lemma \ref{lemme.maj.PPsi})\\
Let $\tilde{P}$ denote the operator $\tilde{P}:\mc{E}^*\longrightarrow\mc{F}'$ defined by 
for all $U\in\mc{E}^*$ and for all $\varphi\in\mc{F}$
\begin{equation}\label{def.Ptilde}
\lang \tilde{P}U,\varphi\rang:=\lang U\gauchedroite{\D_s},\varphi\rang
+\int\dd\tau \lang U(\tau),(\Box+m^2)\varphi(\tau)\rang
\end{equation}
and for $k\in\N^*$ and $I\subset\lent 1,k\rent$ we denote by $\tilde{P}^k_I$ the unique continuous operator 
$\tilde{P}^k_I:\mc{E}^{k*}\longrightarrow\mc{L}_k(\mc{F})$ such that for all $U=U_1\otimes\cdots U_k\in(\mc{E}^*)^{\otimes k}$ and 
for all $\varphi=(\varphi_1,\ldots,\varphi_k)\in\mc{F}^k$
$$
\lang \tilde{P}^k_I U,\varphi\rang=
\prod_{i\in I}\lang \tilde{P} U_i,\varphi_i\rang \prod_{j\not\in I}\int_0^T \lang U_j(\tau_j),\varphi_j(\tau_j)\rang\dd\tau_j
$$
Then since $\forall b\in T(2)$, $b\neq\varepsilon$ and $\forall\alpha\in\{0,1\}^{\NO{b}}$ the operator $\frac{\D^{\alpha}\Psi(b)}{\D t^\alpha}$ satisfies 
$\frac{\D^{\alpha}\Psi(b)}{\D t^\alpha}(t)=0$
for all $t\in\cup_{j=0}^{\NO{b}-1}[0,T]^{j}\times\{0\}\times[0,T]^{\NO{b}-1-j}$ we have $P^{\NO{b}}_I\Psi(b)=\tilde{P}^{\NO{b}}_I\Psi(b)$ for all $b\in T(2)$, 
$b\neq \varepsilon$ and for all $I\subset\lent 1,\NO{b}\rent$. 

Let $\varphi\in\mc{E}$ be a solution of \refeq{E} and $U\in\mc{E}^*$ then $(\Box+m^2)\varphi=-\lambda\varphi^2$ and using property \ref{Hq.algebre} and 
definition \refeq{def.Ptilde} we get 
\begin{equation}\label{base.dem.lemme.PPsi}
\begin{split}
\left\vert 
\lang \tilde{P}U,\varphi\rang
\right\vert
\le&
\NO{U}_{*1}\left(\NO{\varphi(s)}_{H^q}+\NO{\dsurd{\varphi}{t}(s)}_{H^q}\right)+
C_q\vert\lambda\vert T\NO{U}_{*1}\NO{\varphi}_{\mc{E}}^2\\
\le & \left(2+C_q\vert\lambda\vert T\NO{\varphi}_{\mc{E}}\right)\NO{\varphi}_{\mc{E}}\NO{U}_{*1}
\end{split}
\end{equation}

Let $b\in T(2)$, $\NO{b}=:k\ge 2$, then consider a sequence $U_n=(U^{(1)}_n\otimes\cdots\otimes U^{(k)}_n)\in (\mc{E}^*)^{\otimes k}$ such that 
$U_n\to\Psi(b)$ in $\mc{E}^{*k}$. Let $\varphi\in\mc{E}$ be a solution of \refeq{E} and $I\subset\lent 1,k\rent$ then using 
\refeq{base.dem.lemme.PPsi} and property \ref{Hq.algebre} we get 
$$
\left\vert 
\lang \tilde{P}^k_IU_n,\varphi\rang
\right\vert
\le
\left(2+C_q\vert\lambda\vert T\NO{\varphi}_{\mc{E}}\right)^{\vert I\vert}\NO{\varphi}_{\mc{E}}^{\vert I\vert} (T C_q)^{k-\vert I\vert} 
\NO{\varphi}_{\mc{E}}^{2k-\vert I\vert}
\NO{U_n}_{*k}
$$
then taking the limit $n\to\infty$ in the last inequality we get 
$$
\left\vert 
\lang \tilde{P}^k_I\Psi(b),\varphi\rang
\right\vert
\le
\left(2+C_q\vert\lambda\vert T\NO{\varphi}_{\mc{E}}\right)^{\vert I\vert}\NO{\varphi}_{\mc{E}}^{\vert I\vert} (T C_q)^{k-\vert I\vert} 
\NO{\varphi}_{\mc{E}}^{2k-\vert I\vert}
\NO{\Psi(b)}_{*k}
$$
hence lemma \ref{lemme.RdC} completes the proof of lemma \ref{lemme.maj.PPsi}.
\end{proof}

Then using \refeq{base.dem.2}, lemma \ref{lemme.maj.PPsi} and the fact that the number $p_k$ 
of Planar Binary Tree $b$ such that $\vert b\vert=k$ satisfies 
$p_k\le 4^k$
we get 
\begin{multline}\label{base.dem.3}
\vert \Delta_N\vert \le\\
\NO{\psi}_{\mc{E}^*}
\sum_{\beta=N+1}^{2N-1}\vert \lambda\vert^{\beta-1}
\sum_{\substack{1\le l\le k\le N\\ k+l=\beta}}
4^{k-1}
C^{k-l}_k
\NO{\varphi}^{k+l}_{\mc{E}}(C_qT)^{l}(C_qMT)^{k-1}
\left(2+\vert \lambda\vert C_qT\NO{\varphi}_{\mc{E}}\right)^{k-l}
\end{multline}
Let $A$ denotes the quantity $A:=\vert\lambda\vert C_q T\NO{\varphi}_{\mc{E}}$ then using this notation \refeq{base.dem.3} writes
$$
\vert \Delta_N\vert \le
\NO{\psi}_{\mc{E}^*}\NO{\varphi}_{\mc{E}}
\sum_{\substack{1\le l\le k\le N\\ N+1\le k+l\le 2N-1}}
C^{k-l}_k (4MA)^{k-1}A^{l}
\left(2+A\right)^{k-l}
$$
But for all $(k,l)\in(\N^*)^2$ such that $1\le l\le k\le N$ and $N+1\le k+l\le 2N-1$ we have $k\ge [N/2]$, so we have
$$
\vert \Delta_N\vert \le
\NO{\psi}_{\mc{E}^*}\NO{\varphi}_{\mc{E}}
\sum_{k=[N/2]}^N
(4MA)^{k-1}
\left(2+2A\right)^{k-1}
$$
But since \refeq{condition.fin} is satisfied we get $8MA(1+A)<1$ hence the last inequality shows that $\Delta_N$ tends to $0$ when $N$ tends to 
infinity which completes the proof of theorem \ref{prop.fin}.

\appendix
\section{Appendix: Planar Binary Trees}\label{app.1}
Here we will prove the lemma \ref{noeud.dem}. Let begin with the first part of the lemma which is equivalent to
\begin{lemma}\label{combinatoire}\sl{
Let $b$ belong to $T(2)$, $\NO{b}=\beta$ ($\beta\ge 2$), then we have
\begin{equation}\label{somme.combinatoire}
\sum_{\substack{0\le l\le k\le \beta\\k+l=\beta}}\sum_{\substack{a\in
  T(2), \NO{a}=k\\E\in\{\varepsilon,Y\}^k\vert n_Y(E)=l\\\text{such that }E\propto
  a=b}}(-1)^{\vert a\vert}=0
\end{equation}
}
\end{lemma}
\begin{proof}
Let $b\in T(2)$ be such that $\NO{b}=\beta$ ($\beta\ge 2$). Then let us denote by $L$
the integer defined by
\begin{multline*}
L:=\max\{i\in\N\text{ such that }\\ \exists a\in T(2), \NO{a}=\beta-i,\ \exists
E\in\{\varepsilon,Y\}^{\beta-i} \text{ such that }n_Y(E)=i \text{ and } E\propto a=b\}
\end{multline*} 
Then since $\beta\ge 2$ we have $L\ge 1$. 
Define $K$ by $K:=\beta-L$ 
and let $A\in T(2)$, $\NO{A}=K$ and $\hat E\in\{\varepsilon,Y\}^K$ such that $\hat E\propto A=b$ (and then necessarily $n_Y(\hat E)=L$). 
Note that $A$ is actually unique: it is obtained by removing all pairs of $b$ which are sons of the same vertex.
Let denote by $I$ the set of indices $1\le i\le K$ such that $E_i=Y$. Then for all $J\subset I$ we will denote by $E^J$ the $K$-uplet $E^J:=(E^J_1,\ldots,E^J_K)$
where for all $j$ in $\lent 1,K\rent$, $E^J_j$ defined by 
\begin{equation*}
E^J_j=
\begin{cases}
Y&\text{if }j\in J\\
\varepsilon&\text{if }j\in \lent 1,K\rent\setminus J
\end{cases}
\end{equation*}

$\forall j\in J$, $E^J_j=Y$ and $\forall i\in \lent 1, K\rent \setminus J$ $E^J_i=\varepsilon$.

Let $k,l$ be some integers such that $1\le l\le k\le \beta$ and $k+l=\beta$. Then for all $a\in T_k$ such that there exists $E_a\in\{\varepsilon,Y\}^k$
which satisfies $b=E\propto a$, there is an unique subset $J\subset I$ such that $a=E^J\propto A$ and then we have $k\ge \vert J\vert =l\ge 1$.
In the other hand for all $J\subset I$ such that $\vert J\vert \ge 1$ there exists an unique $\tilde E\in\{\varepsilon,Y\}^{K+\vert J\vert}$, $n_Y(\tilde E)\ge 1$ 
such that $\tilde E\propto(E^J\propto A)=b$. Hence we have 
$$
-\sum_{\substack{0\le l\le k\le \beta\\k+l=\beta}}\sum_{\substack{a\in
  T(2), \NO{a}=k\\E\in\{\varepsilon,Y\}^k\vert n_Y(E)=l\\\text{such that }E\propto
  a=b}}(-1)^{\vert a\vert}=0
 -\sum_{\substack{J\subset I\\ \vert J\vert\le L}}(-1)^{\vert E^J\propto A\vert}
$$
but $\vert E^J\propto A\vert=K+\vert J\vert -1$ so the previous equality leads to
\begin{equation*}
-\sum_{\substack{0\le l\le k\le \beta\\k+l=\beta}}\sum_{\substack{a\in
  T(2), \NO{a}=k\\E\in\{\varepsilon,Y\}^k\vert n_Y(E)=l\\\text{such that }E\propto
  a=b}}(-1)^{\vert a\vert}=-\sum_{l=0}^{L}C^l_L(-1)^{K+l-1}
  =
  (-1)^{K}(1-1)^L=0
\end{equation*}
which completes the proof. 
\end{proof}

Let focus on the second part of lemma \ref{noeud.dem}
\begin{lemma}\label{combinatoire2}\sl{
Let $p\in\N^*$, $a\in T(2)$ be such that $\NO{a}=p$ and $E\in\{\varepsilon, Y\}^p$ then we have $\NO{E\propto a}=p+n_Y(E)$ and
\begin{equation}\label{ce.quil.faut.avoir}
\lang P^{p+n_Y(E)}_{\lent 1,p+n_Y(E)\rent}\Psi(E\propto a), (\varphi,\ldots,\varphi)\rang=
\lang P^p_{I_E}\Psi(a),(\varphi^{\alpha^{I_E}_1},\ldots,\varphi^{\alpha^{I_E}_p}
\rang
\end{equation}
where $I_E:=\{j\in\lent 1,p\rent\text{ such that }E_j=\varepsilon\}$ and $\alpha^{ I_E}_j=2$ if $j\not\in I_E$ and $1$ otherwise. 
}
\end{lemma}
\begin{proof}
Let $k\in\N^*$ and $U$ belong to $\mc{E}^{*k}$, then for all $K\subset\lent 1,k\rent$, for all $t^{\vee K}\in[0,T]^{k-\vert K\vert}$ and 
for all $f^{\vee K}\in (H^q)^{k-\vert K\vert}$ we consider the element $U^{\vee K}(t^{\vee K},f^{\vee K})$ of 
$\mc{E}^{*\vert K\vert}$ defined by $\forall \tau\in[0,T]^{\vert K\vert}$ and $\forall g\in (H^q)^{\vert K\vert}$, 
$
\lang U^{\vee K}(t^{\vee K},f^{\vee K})(\tau),g\rang:=\lang U(\tilde{t}),\tilde{f}\rang
$ 
where $\tilde{t}$ and $\tilde{f}$ are defined by 
\begin{equation}\label{def.Psi^j}
\left\{
\begin{array}{l}
\tilde{t}_r:=t^{\vee K}_{v(r)}\text{ if }r\not\in K \\
\tilde{t}_r:=\tau_{k(r)}\text{ if }r\in K
\end{array}
\right.
\text{ and }
\left\{
\begin{array}{l}
\tilde{f}_r:=f^{\vee K}_{v(r)}\text{ if }r\not\in K \\
\tilde{f}_r:=g_{k(r)}\text{ if }r\in K
\end{array}
\right.
\end{equation}
here $v(r):=\text{card} \{k\le r\text{ such that }k\not \in k\}$ and 
$k(r):=\text{card} \{k\le r\text{ such that }k \in K\}$.

First we will treat the case $n_Y(E)=1$ then we will see how to generalize the result.
For $j\in\lent 1,k\rent$ we define $E^{(j,k)}=(E^{(j,k)}_1,\ldots,E^{(j,k)}_k)\in\{\varepsilon,Y\}^k$ by 
$E^{(j,k)}_r=\varepsilon$ if $r\neq j$ and $E^{(j,k)}_j=Y$. Let $t\in[0,T]^{k-1}$ and $(f_1,\ldots,f_{k-1})\in (H^q)^k$ then we consider the 
element $\Psi(a)^{\vee \{j\}}(t,f)$ of $\mc{E}^{*}$. 
In view of the definition of $\Psi(b)=\DD(b)\psi$ we have $\Psi(E^{(j,k)}\propto a)=\DD\left[\Psi(a)^{\vee \{j\}}(t,f)\right]\in \mc{E}^{*2}$. 
Then the calculations done in the proof of proposition \ref{prop.ordre2} shows that for all $\varphi\in\mc{E}$ solution of \refeq{E}  we have
\begin{equation}\label{cas.nY=1}
\lang P^2_{\lent 1,2\rent}\DD\left[\Psi(a)^{\vee \{j\}}(t,f)\right],(\varphi,\varphi)\rang
=
\int_0^T
\lang \Psi(a)^{\vee j}(t,f)(\tau),\varphi^2(\tau)\rang\dd\tau
\end{equation}
Hence using the definition \refeq{def.Psi^j} of $\Psi(a)^{\vee \{j\}}(t,f)$ we find that the lemma is true if $E=E^{(j,k)}$ \ie when $n_Y(E)=1$.

Let $M\in\N^*$ 
and $E\in\{\varepsilon,Y\}^k$ be 
such that $n_Y(E)=M$. Then we define $J_E\subset\lent 1,k\rent $ as the set of indices $j\in\lent 1,k\rent$ such that $E_j=Y$, then 
since $n_Y(E)=M$ we have $\vert J_E\vert=M$. We denote $J_E:=\{j_1,\ldots,j_M\}$ where $j_M<j_{M-1}<\cdots <j_1$. 
Then one can show easily that we have 
$$
b:=E\propto a=E^{(j_M,k+M-1)}\propto(E^{(j_{M-1},k+M-2)}\propto(\cdots\propto(E^{(j_1,k)}\propto a))\cdots)
$$
Hence if we denote by $a_1$ the Planar Binary Tree $a_1:=E^{(j_{M-1},k+M-2)}\propto(\cdots\propto(E^{(j_1,k)}\propto a))\cdots)$ we have 
$b=E\propto a=E^{(j_M,k+M-1)}\propto a_1$. Then for all  
\begin{equation*}
\begin{split}
t^{\{j_M\}}&=(t_1,\ldots,t_{j_M-1},t_{j_M+1},\ldots,t_{k+M-1})\in[0,T]^{k+M-2}\\
f^{\{j_M\}}&=(f_1,\ldots,f_{j_M-1}, f_{j_M+1},\ldots,f_{k+M-1})\in (H^q)^{k+M-2}
\end{split}
\end{equation*}
we can use \refeq{cas.nY=1} and the fact that 
$$
\DD\left[\Psi(a_1)^{\vee j_M}(t^{\{j_M\}},f^{\{j_M\}})\right]=\Psi(b)^{\vee \{j_M,j_M+1\}}(t^{\{j_M\}},f^{\{j_M\}})
$$ 
in order to obtain 
\begin{multline*}
\lang P^2_{\lent 1,2\rent}\Psi(b)^{\vee \{j_M,j_M+1\}}(t^{\{j_M\}},f^{\{j_M\}}),(\varphi,\varphi)\rang
=\\
\int_0^T
\lang \Psi(a_1)(t),
(f_1,\ldots,f_{j_M-1},\varphi^2(t_{j_M}),f_{j_M+1},\ldots,f_{k+M-1})\rang\dd t_{j_M}
\end{multline*}
where $t$ denotes the $(k+M-1)$--uplet $t:=(t_1,\ldots,t_{j_M-1},t_{j_M},t_{j_M+1},\ldots,t_{k+M-1})$. Then writing 
$a_1=E^{(j_{M-1},k+M-2)}\propto a_2$ one can use the same arguments to show that 
\begin{multline*}
\lang P^4_{\lent 1,4\rent}\Psi(b)^{\vee \{j_M,j_M+1,j_{M-1}+1,j_{M-1}+2\}}(t^{\{j_M,j_{M-1}\}},f^{\{j_M,j_{M-1}\}}),(\varphi,\varphi,\varphi,\varphi)\rang
=\\
\iint_{[0,T]}\dd t_{j_M}\dd t_{j_{M-1}}
\lang \Psi(a_2)(t),
(f_1,\ldots,\varphi^2(t_{j_M}),f_{j_M+1},\ldots,f_{j_{M-1}-1},\varphi^2(t_{j_{M-1}}),\ldots,f_{k+M-1})\rang
\end{multline*}
Hence 
Doing this operation successively for $j_{M-2},\ldots,j_1$ we finally get 
\begin{equation}\label{fin.de.fin}
\lang P^{\vert K\vert}_{\lent 1,\vert K\vert\rent}\Psi(b)^{\vee K}(t^{\vee K},f^{\vee K}),(\varphi,\ldots,\varphi)\rang
=
\iint_{[0,T]^{M}}\dd t_{j_M}\cdots\dd t_{j_1}
\lang \Psi(a)(t),
(\tilde{g}_1,\cdots,\tilde{g}_k)\rang
\end{equation}
where $K:=\bigcup_{r=1}^M\{j_r+M-r,j_r+M-r+1\}$ and where $\tilde{g}_r:=\varphi^2(t_r)$ if $r\in J_E$ and 
$\tilde{g}_r:=f_{v(r)}$ otherwise. Hence, considering the element 
$\lang P^{\vert K\vert}_{\lent 1,\vert K\vert\rent}\Psi(b)^{\vee K}(\cdot,\cdot),(\varphi,\ldots,\varphi)\rang$ of 
$\mc{E}^{*(k+M-2M)}$ and using \refeq{fin.de.fin}, we finally get 
\begin{equation*}
\lang P^{k+M}_{\lent 1,k+M\rent}\Psi(b),(\varphi,\ldots,\varphi)\rang
=
\lang P^{k}_{\lent 1,k\rent\setminus J_E}\Psi(a),
(\tilde{h}_1,\cdots,\tilde{h}_k)\rang
\end{equation*}
where $\tilde{g}_r:=\varphi^2$ if $r\in J_E$ and $\tilde{h}_r:=\varphi$ otherwise \ie we obtain exactly identity 
\refeq{ce.quil.faut.avoir}.
\end{proof}



\end{document}